\documentclass[11pt, a4paper]{amsart}
\usepackage[utf8]{inputenc}
\usepackage[english]{babel}
\usepackage{graphicx} 
\usepackage{amsmath,amsfonts,amssymb,mathtools}
\usepackage{array}
\usepackage{hyperref}
\usepackage{url}
\usepackage{caption}
\usepackage{xcolor}
\newtheorem{thm}{Theorem}[section]
\newtheorem{prop}[thm]{Proposition}

\newtheorem{rem}[thm]{Remark}

\newcommand{\zek}{Zekovi\'{c}}
\newcommand{\hhh}{\mathrm{H}(2)}

\title[An exploration of low-crossing and chiral cosmetic bands]{An exploration of low crossing and chiral cosmetic bands with grid diagrams}

\author{Agnese Barbensi, Daniele Celoria, Christopher Ktenidis \\with an appendix by Kazuhiro Ichihara, In Dae Jong, Masakazu Teragaito}
\address{School of Mathematics and Physics, the University of Queensland, 4072 Brisbane Australia.}
\email{a.barbensi@uq.edu.au}
\email{d.celoria@uq.edu.au}
\email{c.ktenidis@student.uq.edu.au}

\address{Department of Mathematics, College of Humanities and Sciences, Nihon University, 3-25-40 Sakurajosui, Setagaya-ku, Tokyo 156-8550, JAPAN}
\email{ichihara.kazuhiro@nihon-u.ac.jp}

\address{Department of Mathematics, Kindai University, 3-4-1 Kowakae, Higashiosaka City, Osaka 577-0818, Japan} 
\email{jong@math.kindai.ac.jp}

\address{Department of Mathematics and Mathematics Education, Hiroshima University, 1-1-1 Kagamiyama, Higashi-hiroshima 7398524, Japan.}
\email{teragai@hiroshima-u.ac.jp}

\date{}

\begin{document}

\begin{abstract}
We computationally explore non-coherent band attachments between low crossing number knots, using grid diagrams. We significantly improve the current $\hhh$-distance table. In particular, we find two new distance one pairs with at most seven crossings: one between $3_1\#3_1$ and $7_4m$, and a chirally cosmetic one for $7_3$. We further determine a total of 33 previously unknown $\hhh$-distance one pairs for knots with up to $8$ crossings. The appendix by Kazuhiro Ichihara, In Dae Jong and Masakazu Teragaito contains a construction explaining the existence of chirally cosmetic bands for an infinite family of knots, including $5_1,\, 7_3$ and $8_8$.
\end{abstract}

\maketitle
% \keywordss{Knot theory, Band attachments, $H(2)$-moves.\\ \textbf{Corresponding author:} Daniele Celoria, University of Queensland,\\ \texttt{d.celoria@uq.edu.au}}
\section{Introduction}\label{sec:introduction}
Given a knot $K \subset S^3$, one can obtain a new one by a local move called \textit{non-coherent band attachment}~\cite{recent_vaz,kane3} or \textit{$\hhh$-move}~\cite{hoste}. Similar to crossing changes, these moves are unknotting operations and therefore induce a well-defined metric on knots, referred to as $\hhh$-distance~\cite{kanenobu}. 

The first effort to systematically study $\hhh$-moves is due to Lickorish, who produced the first table of knots with $\hhh$-distance $1$ from the unknot~\cite{lickorish}. Since then, several methods to obstruct low distance have been found (see~\textit{e.g.}~\cite{kane1, kane2, kane3, kai, recent_vaz}). Over the past decade, non-coherent band attachments have gained significant attention, partially due to their relevance in biological applications~\cite{recent_vaz}. 

In this paper, we perform a computational exploration of $\hhh$-moves for low crossing number knots, using grid diagrams. Grid diagrams are a versatile discrete analogue of knot diagrams. They are widely used throughout low-dimensional topology and contact geometry~\cite{cromwellI,NGgrids}, and are the starting point to define a combinatorial counterpart to knot Floer homology~\cite{SOS}. 

A deep result by Dynnikov~\cite{dinnikov} implies that grid diagrams for the unknot are theoretically guaranteed to simplify \textit{monotonically}, \textit{i.e.}~without increasing complexity. While this result fails to extend to other knot types, in practice, the same method applies in general: performing local simplification moves on a grid is overwhelmingly likely to result in a diagram minimising the grid size within its knot type. This experimentally observed behaviour makes knot recognition more manageable on average than in other computational knot models. Furthermore, grid diagrams are extremely easily randomised. These two characteristics make them a computationally convenient model to investigate the behaviour of random knots and operations on them~\cite{griglie_dna, gridpym}. ~\\

Here, we leverage on these features to look for $\hhh$-moves between knots with small crossing number. Specifically, we focus on knots with up to $8$ crossings; we generate random diagrams for each of these knots, perform all available $\hhh$-moves on the resulting configurations, and identify the knots obtained. As a result, we find previously unknown bands between pairs of knots, thus improving the current $\hhh$-distances table. The state of the art for such a table for knots with at most $7$ crossings was compiled by Kanenobu~\cite{kane3}. Therein, new results by \zek~\cite{zekovic} are included, such as single $\hhh$ moves connecting $7_1$ to $5_2m, 7_5m, 7_7m$, and between $6_2$ and $7_2$. Here, and in the rest of the paper, the $m$ denotes mirrors, consistently with KnotInfo's notation~\cite{knotinfo}. 

We find a total of $322$ distance one pairs (including the $43$ cosmetic ones). These include two previously unknown between knots with fewer than $7$ crossings: $7_4 \leftrightarrow 3_1m\#3_1m$, for which it was previously only known that the distance was between 1 and 3, and $7_3 \leftrightarrow 7_3m$. Some --but not all-- of the ones involving knots with $8$ crossings are mentioned in \zek's PhD thesis~\cite{zekovicPHD}, but different conventions and language barriers made it hard to properly compare our datasets. In our computations, we find a total of $33$ pairs that were not found by \zek. ~\\

Arguably, \zek's main result in~\cite{zekovic}  was the discovery of a \textit{chirally cosmetic} band attachment for $5_1$. This is a $\hhh$ move taking one knot type into its mirror. Such band attachments are of interest due to their connections with chirally cosmetic surgeries, via the Montesinos trick~\cite{montesinos}. To date, the only known chirally cosmetic crossings on knots with fewer than $9$ crossings are the aforementioned $5_1$ (\zek~\cite{zekovic}) and $8_8, 8_{20}$ (Moore-Vazquez~\cite{recent_vaz}). As mentioned above, here we find a single $\hhh$ move connecting $7_3$ to its mirror.

Finally, it was noted in~\cite{recent_vaz} that the occurrence of chiral cosmetic band attachments is extremely low, as a fraction of the sampled population. While this is still the case for our grid-based model, we note that the relative probabilities are substantially higher.\\

The appendix to this paper contains a way to illustrate the chirally cosmetic banding on $7_3$ described above. Furthermore, it  generalises this construction to an infinite family of knots --including, in particular, the knots $8_8$ and $5_1$-- providing a broader framework for understanding chirally cosmetic bandings in these cases.

\subsection*{Acknowledgements} The authors are grateful for the funding provided by UQ's School of Mathematics and Physics and by SEES (Student Enrichment \& Success) as part of the Winter Research Project that kickstarted this paper. We also wish to thank Ben Burton, Craig Hodgson, Neil Hoffman for suggestions and interesting discussions.

\section{Grid diagrams and band attachments}\label{sec:grids}

\begin{figure}
\centering
\includegraphics[width=0.5\linewidth]{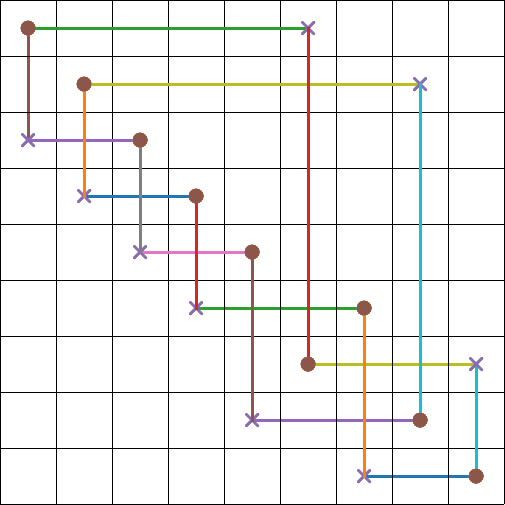}
\caption{A grid diagram of size $9$, representing the knot $7_3$. }
\label{fig:7_3grid}
\end{figure}

A \textit{grid diagram of size }$n$ is an $n\times n$ square grid, divided into $n^2$ small squares. Two sets of $n$ \textit{markings}, usually denoted by $X$ and $O$, are placed within the grid, so that no two markings of the same kind are in the same row or column, and each small square contains at most one marking (see Figure~\ref{fig:7_3grid}). 
A knot diagram is then obtained by connecting the markings on the same row and column, and resolving each crossing as an overpass for the vertical strand. The size of the grid is also referred to as the grid number. An orientation of the link is encoded by following the convention dictating that the arcs go from $O$ to $X$ markings horizontally, see \textit{e.g.}~\cite{SOS,gridpym} for more details. Just as for Reidemeister moves and regular diagrams, there is a finite set of \textit{grid moves} between two grid diagrams if and only if they represent the same knot type. These moves are known collectively as \textit{Cromwell moves}.\\

\begin{figure}[h!]
\captionsetup{width=.85\linewidth}
\centering
\includegraphics[width=1.1\linewidth]{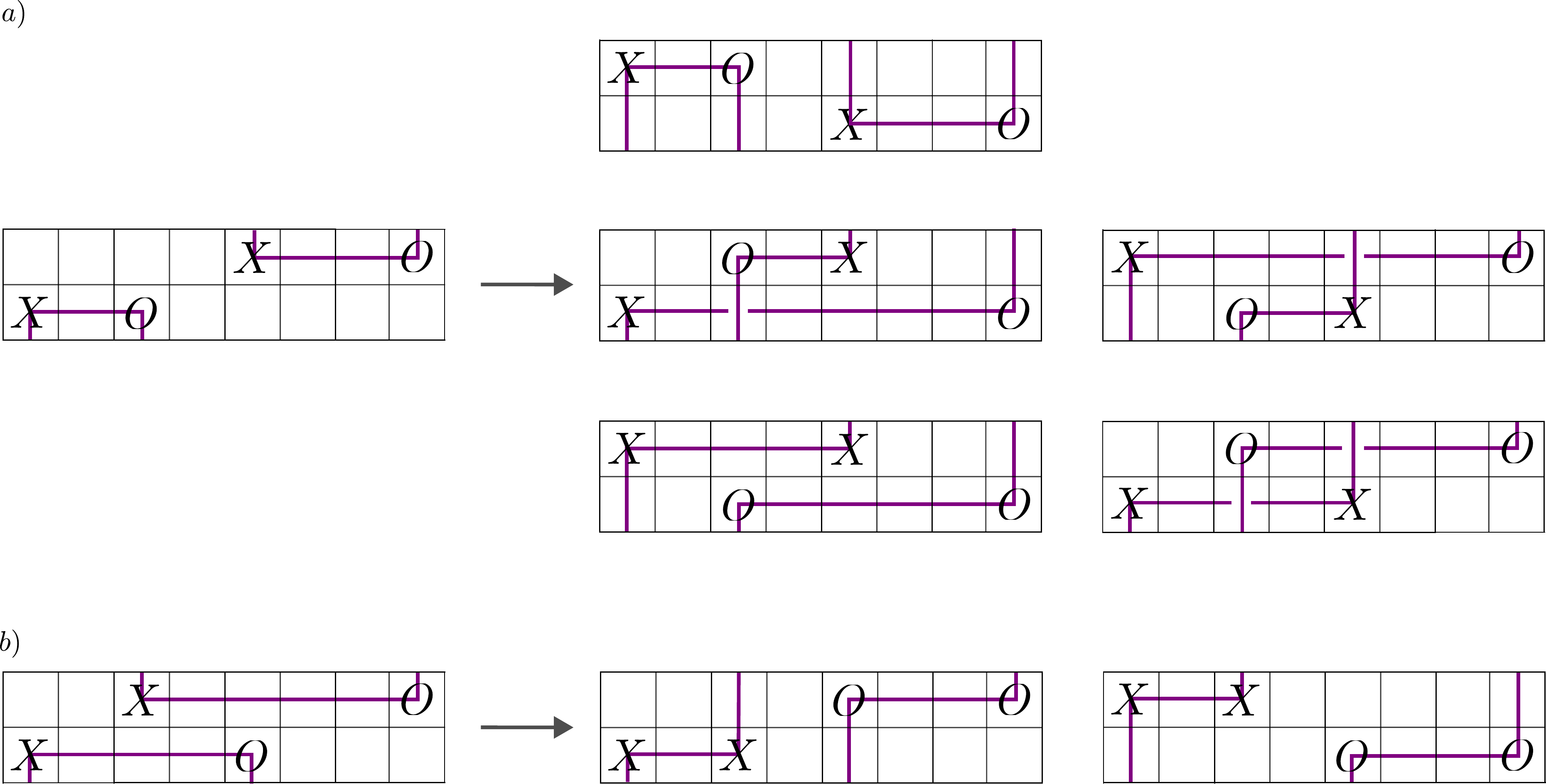}
\caption{\textit{(a)} Depicts a single portion of a grid diagram (on the left) being changed in
each of the five possible ways (on the right). The top row shows that the link diagram on the left changes by a planar isotopy if both pairs of markings are moved. The mid row shows that swapping markings of the same type in the left diagram results in a coherent band attachment. Note that the two possibilities produce equivalent diagrams. The bottom row shows that swapping two markings of different types in the left diagram  results in a non-coherent band attachment. Note that in this case the two possibilities might produce non-equivalent links. \textit{(b)} The two moves on the left diagram correspond to non-coherent band attachments yielding equivalent links.}
\label{fig:band_grids}
\end{figure}

A \textit{band attachment} (also known as band surgery or $\hhh$ move) on a link $L \subset S^3$ is the following operation: given a smooth embedding of the square $\iota\colon [0,1]^2 \hookrightarrow S^3$ such that $L \cap \iota([0,1]^2) = \iota([0,1]\times \{0,1\})$, a new link $L_\iota$ is obtained as  
\begin{equation*}
L_\iota = (L \setminus \iota ([0,1]\times \{0,1\})) \cup \iota (\{0,1\} \times [0,1]).
\end{equation*}
A band attachment is called \textit{coherent} if the two links it relates have different numbers of  components, and \textit{non-coherent} otherwise.

There are combinatorial local moves that carry out (non-)coherent band attachments on grid diagrams. While these are known to experts (see \textit{e.g.}~\cite{SOS}), we include an explicit description here, as they can't be readily found in the literature.

Consider two consecutive rows in a grid (the case for columns works in the same way);
each contains one $X$ and one $O$ marking. We will assume for simplicity (but without loss of generality) that the markings on the lower row are not adjacent to any marking on the upper row.

We say that two consecutive rows are \textit{interleaved} if the arcs spanned by the markings partially overlap (see the left diagram in Figure~\ref{fig:band_grids}(b)), and \textit{non-interleaved} otherwise.

We can non-trivially move pairs of markings vertically within the two rows in five possible ways, while still producing a valid grid diagram. The first possibility is to move one connected XO pair of markings upwards and the other pair downwards. The effect on the link is either a crossing change or equivalent to a planar isotopy on the underlying knot diagram. These two alternatives depend on whether the two markings were interleaved or not, see the top row of Figure~\ref{fig:band_grids}(a) for the non-interleaving case. The second possibility is to move a single marking upward and another non-connected single marking
downward and leave the other two in the same position. In this case the effect of the move depends on whether we move two markings of the same kind (\textit{e.g.}~two $O$s), or not. In the first case, the move corresponds to a coherent band attachment, and in the second case to a non-coherent one, see the mid and bottom row in Figure~\ref{fig:band_grids}(a). Note that in the latter case the result is an invalid grid $G$ as the link orientation is not preserved by the move. However, it is possible to exhibit a valid grid $G'$ whose unlabelled markings are the same as $G$. In our computations, non-coherent band attachments are achieved using built-in functions in \texttt{GridPym}, see~\cite{gridpym} and Section~\ref{sec:data_availability} for further details and documentation. Depending on whether the rows/columns are interleaved or not, and whether the corresponding arcs are oriented in the same direction or not, exchanging one pair $X \leftrightarrow O$ among the two rows can yield non-equivalent band attachments, see Figure~\ref{fig:band_grids}(a) mid and bottom row, and Figure~\ref{fig:band_grids}(b) for an example. Note that this can occur in the coherent and non-coherent case alike.

\section{Methods and results}\label{sec:code}
Our computational exploration was achieved by using software in~\cite{repo}, based on GridPyM~\cite{gridpym}, see Software and data availability (Section~\ref{sec:data_availability}). \\Here we focus on the $43$ knots with crossing numbers $\le 8$ (up to mirroring, and including connected sums). Note that all these knots have grid number $\le 10$.

\begin{figure}[h!]
\centering
\includegraphics[width=1\linewidth]{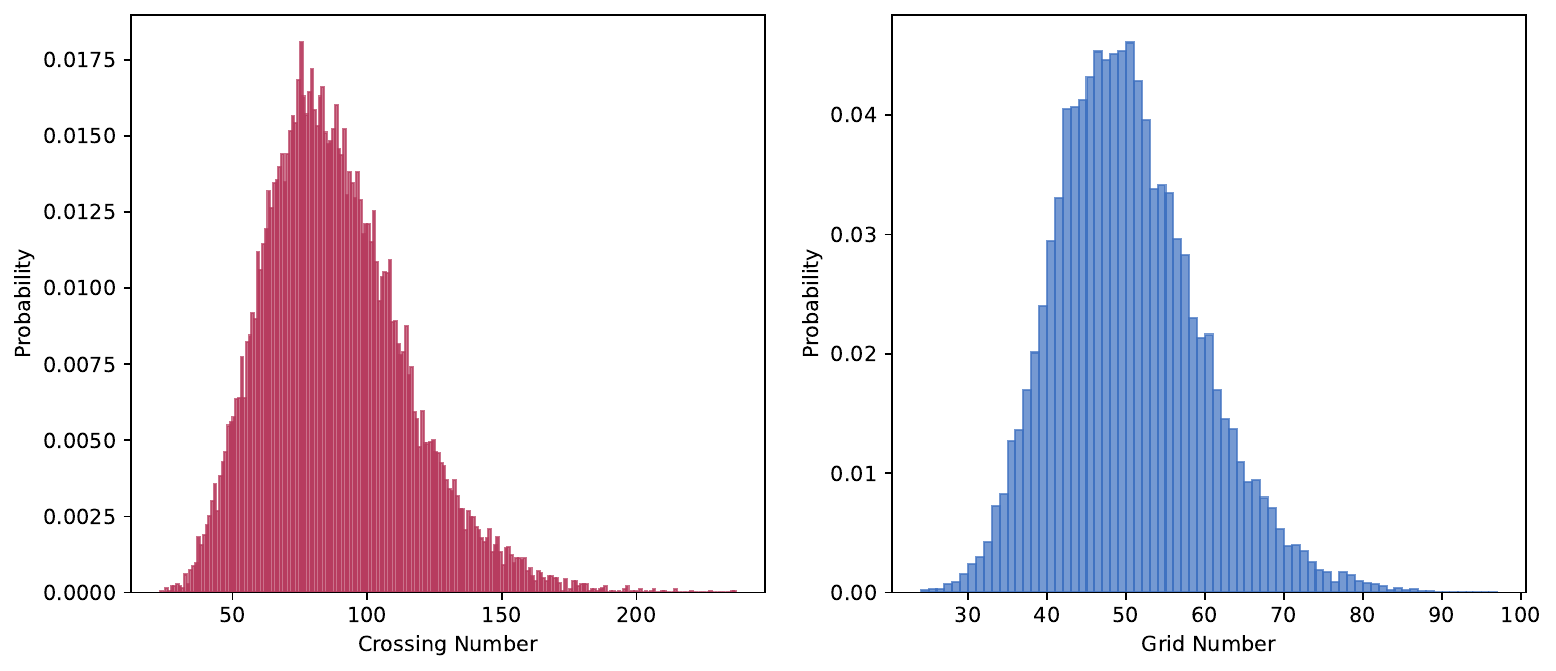}
\caption{Distribution of crossing number and grid number for grids generated by a standard randomisation run: $1000$ random moves and $400$ distinct randomisations per seed diagram.}
\label{fig:stats_on_grids}
\end{figure}

For every chiral representative (up to orientation) of the knot types considered, we select grid diagrams minimising the grid number. We then randomise each of these seed configurations by applying a customisable amount of random grid moves. This is achieved by using the \texttt{scramble\_grid} function on GridPyM. This function takes as input a grid diagram $G$ and a positive integer $m$. The output is a new grid diagram, representing the same knot type as $G$, obtained by applying a sequence of $m$ random Cromwell moves. For a detailed description of the function and for GridPym's documentation see~\cite{gridpym}.

In a standard run, for each seed configuration $G$, we compute $400$ different randomisations by applying the \texttt{scramble\_grid} function with $m = 1000$. This yields a collection of $400$ complex and high-grid number diagrams, see Figure~\ref{fig:stats_on_grids}. 
%In a standard run, \comm{we apply $m1000$ }different moves a total of $400$ times for each seed configuration, yielding a collection of complex and high-grid number diagrams, see Figure~\ref{fig:stats_on_grids}. On average, this results in performing over $32,000$ band attachments per seed configuration, per run. The number of runs we performed varied depending on the knot type.
Once this collection is created, we search for all the possible non-coherent band attachments in each of the $400$ generated diagrams. We do this by following the built-in methodology explained in Section~\ref{sec:grids}, using GridPym's functions, see Section~\ref{sec:data_availability}, and~\cite{repo, gridpym} for a complete description and further documentation. On average, this results in performing over $32,000$ band attachments per seed configuration, per run. The number of runs we performed varied depending on the knot type.

%Once this collection is created, we perform all possible non-coherent band attachments in each generated diagram, using the built-in methodology explained in Section~\ref{sec:grids}. \comm{Note that on average, this results in performing over $32,000$ band attachments per seed configuration, per run. The number of runs we performed varied depending on the knot type. }

Sage's \texttt{get\_knotinfo()} functionality~\cite{sage} is then used to determine the resulting knot types (following KnotInfo's conventions~\cite{knotinfo}). More details are available in the code repository~\cite{repo}, see also Section~\ref{sec:data_availability}.~\\

While it is often difficult to effectively simplify classical knot diagrams (see~\textit{e.g.}~\cite{hass}), this task is often easier in the grid diagrams model. The software GridPyM implements a simplification function that applies random row/column commutations, correspoonding to planar isotopies or Reidemeister moves, and, whenever possible \textit{destabilisation} moves that decrease the grid's size~\cite{gridpym}. Although monotonic simplification is not guaranteed for non-trivial knots, in practice, this method is likely to result in a minimal grid diagram. Knot recognition is performed after multiple rounds of simplification, thus substantially decreasing computational time. 

\begin{figure}[h!]
\centering
\includegraphics[width=0.3\linewidth]{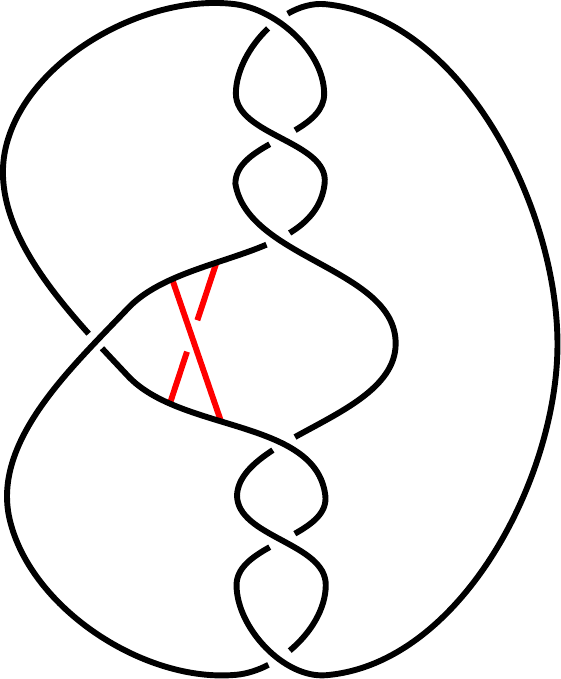}\\~\\
\includegraphics[width=5cm]{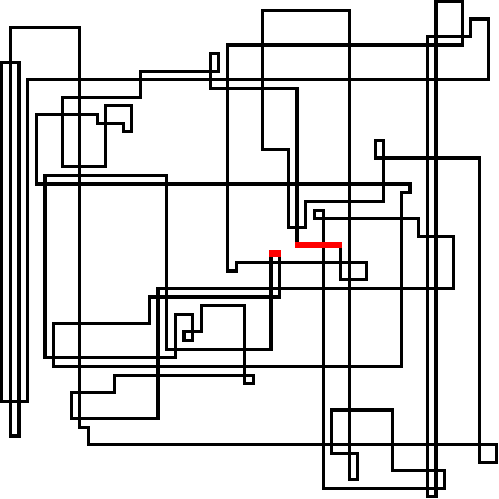}\,\,
\includegraphics[width=5cm]{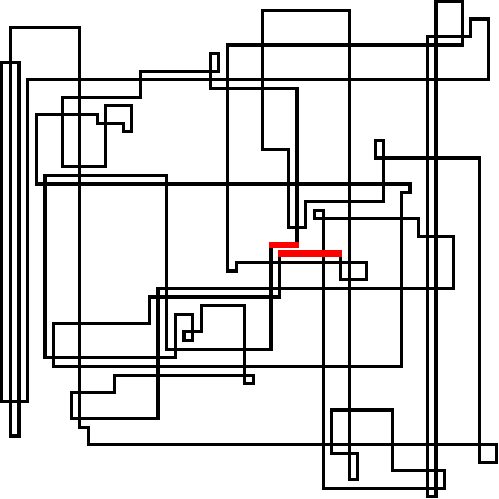}
\caption{\textit{Top:} A non-coherent band attachment from $7_4$ to $3_1m\# 3_1m$. \textit{Bottom:} an explicit example of the band attachment realised in the grids model found by our computations. The band (in red) is attached on row 28.}
\label{fig:newband7_4}
\end{figure}

The combinatorial nature of grid diagrams, along with their randomisation and simplification properties, makes them a simple yet effective tool for computational explorations. Note that all of our computations were performed on standard laptops and are easily parallelisable, thus making similar analysis or extensions fully accessible.\\

\subsection{The $\hhh$-adjacency tables}
We collect a summary of our results in three tables: Tables~\ref{tab:distance_one1}, \ref{tab:distance_one2}, \ref{tab:distance_one3}, containing prime knots with $\le 7$ crossings, prime knots with $8$ crossings and composite knots with $\le 8$ crossings respectively. Overall, we find a total of 33 previously unknown distance one pairs of knots with fewer than $9$ crossings. Two of these pairs connect knots with at most $7$ crossings: $7_4$ with $3_1m \# 3_1m$ (see Figure~\ref{fig:newband7_4}), and $7_3$ with $7_3m$. Note that we are disregarding orientation in our computations. \\

We additionally checked if the data collected could be used to find unknown pairs of knots with $\le 7$ crossings and distance $2$, but didn't succeed. Thus, currently it is still unknown whether the following pairs are at distance $2$ or $3$:
$$(4_1,5_1), (6_1,3_1\#3_1), (6_1 m,3_1\#3_1), (6_3,7_3), (3_1\#3_1, 3_1m\#4_1).$$

Examples of specific band attachments relating each distance one pair (excluding cosmetic ones) are available online~\cite{repo}, see the Software and data
availability, Section~\ref{sec:data_availability}.

\section{Chirally cosmetic band attachments}\label{sec:cosmetic}

Recall that a single non-coherent band attachment converting a knot into its mirror is called \textit{chirally cosmetic}. Such band attachments are related to important open questions in low-dimensional topology, such as the cosmetic surgery conjecture~\cite{ichihara,livingston}. An application of Montesinos' trick (see \textit{e.g.}~\cite[Sec.~6D]{rolfsen}) implies that if two knots are related by a single $\hhh$ move, then their double branched covers are related by a distance one surgery. In the chirally cosmetic case, the surgery connects the filled manifold with its mirror. As mentioned in the introduction,  \zek~exhibited a chirally cosmetic band attachment for $5_1$. Note that in~\cite{noteband_vaz} Moore-Vazquez constructed obstructions to the existence of chirally cosmetic bands for alternating torus knots $T(2,n)$ (with $n$ square-free), using a very elegant argument involving the mapping cone formula in HFK (see also~\cite{lidman}). The square-free condition was later dropped by Livingston~\cite{livingston}, except for $n=9$ (see also~\cite{wu_L9}). 
%Finally, Yang~\cite{yang} and then Ito~\cite{ito} eliminated the $n =9$ case, implying that $5_1$ is the only possible alternating torus knot admitting such a band attachment.

\begin{figure}[h!]
\centering
\includegraphics[width=0.5\linewidth]{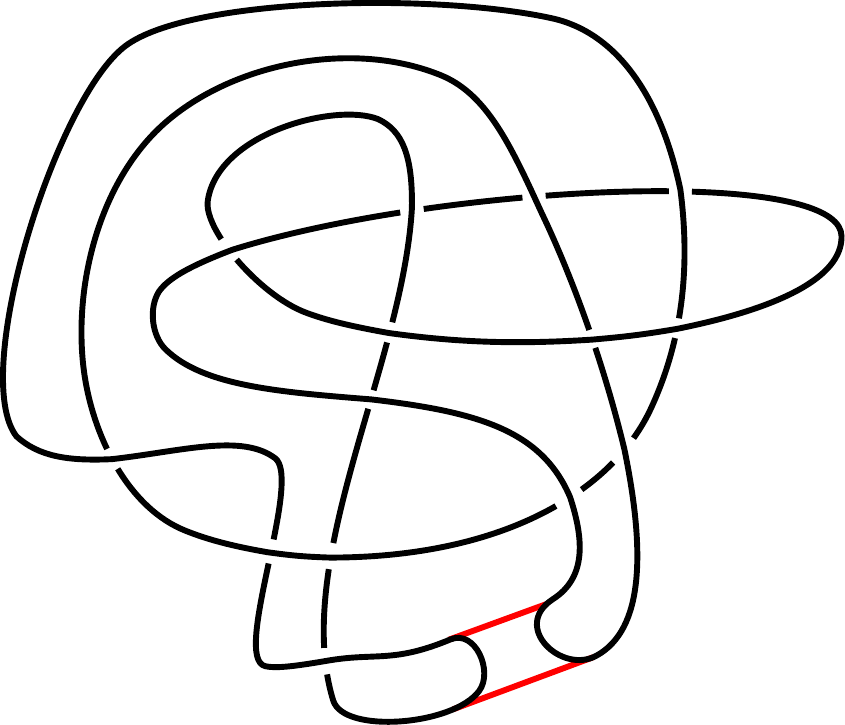}
\caption{A chiral cosmetic band attachment on $7_3m$}
\label{fig:newband7_3}
\end{figure}

Previously, for knots up to $8$ crossing, besides \zek's $5_1$ example, only $8_8$ and $8_{20}$ were known~\cite{recent_vaz} to admit cosmetic smoothings. Using our code, we show the following. 

\begin{prop}\label{prop:chiral}
The knot $7_3$ admits a chirally cosmetic band attachment, displayed in Figure~\ref{fig:newband7_3}.
\end{prop}

Peculiarly, all these small examples (including ours) have been produced through different flavours of computational explorations; so far theory has only been used to obstruct these low-crossing examples, while producing infinitely many examples with at least $9$ crossings~\cite{ichihara}.

\begin{rem}
A result by Livingston~\cite[Thm.~1]{livingston}, together with Proposition~\ref{prop:chiral}, implies that $7_3\#7_3$ is the boundary of a properly embedded M\"obius band in $D^4$.
\end{rem}

The double branched cover of $7_3$ is the lens space $L(13,3)$. By following the recipe beautifully explained in~\cite{owens}, it is possible to explicitly exhibit the knot in $L(13,3)$ providing the lift of the core of the band shown in Figure~\ref{fig:newband7_3}. What follows is a brief description of the procedure. \\

\begin{figure}
\centering
\includegraphics[width=0.8\linewidth]{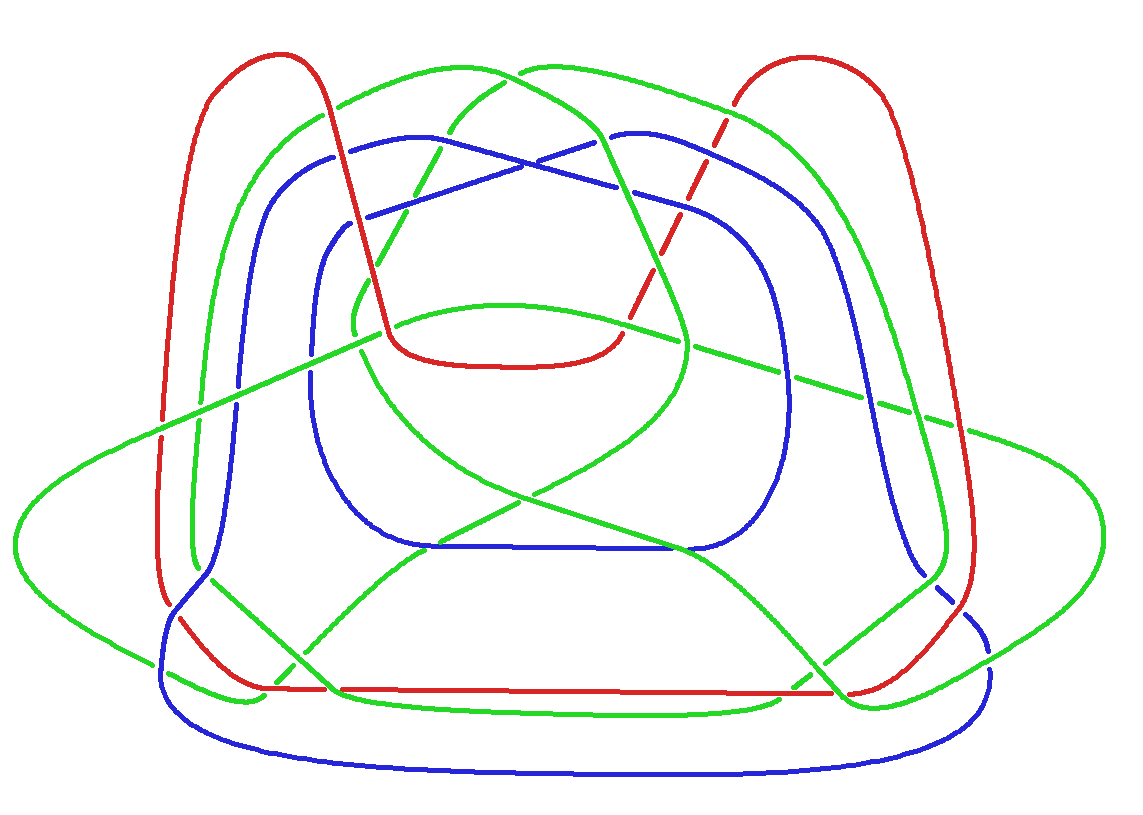}
\caption{A mixed diagram for the knot complement obtained by lifting the band's core to $L(13,3)$. The knot is the green curve, while both the blue and red curves have been filled with framing $\frac{7}{2}$.}
\label{fig:big_mess}
\end{figure}

The unknotting number of $7_3$ is $2$; therefore, we can turn it into the unknot by switching two suitable crossings. This can be done \textit{e.g.}~by performing $(+1)$-framed surgeries over two unknots linking these crossings. We then lift the core of the band and the surgery information to the double branched cover. This produces the mixed surgery diagram shown in Figure~\ref{fig:big_mess}. Using SnapPy~\cite{snappy}, we can reliably identify this as the census manifold \texttt{s463}. Filling the cusp with slopes $(0,1)$ and $(1,0)$, produce $L(13,3)$ and $L(13,10)$ respectively. Note that these slopes are related by an orientation reversing isometry of \texttt{s463}. Additionally, one can explicitly identify the knot in $L(13,3)$ to be (the mirror of) $(13,4)$-surgery on the unknotted component of $L11n219$.\\

To summarise, using \cite[Thm.~4.4]{noteband_vaz} and \cite[Thm.~3]{livingston}, the existence of chirally cosmetic band attachments (for knots with up to $8$ crossings) is still unknown for 
$$6_1, 3_1\#3_1, 7_5, 8_1, 8_2, 8_{13}.$$

We conclude with a table (Table~\ref{tab:chiral_probs}) of the relative occurrences of chiral cosmetic band attachments in the grid model. It is noteworthy to point out that the frequency of occurences is between $1$ (for $8_{20}$) and $3$ (for $8_8$) orders of magnitude higher than in the analogous analysis by Moore and Vazquez~\cite{recent_vaz}. Note that this phenomenon is coherent with the fact that the grid diagram model is \textit{overly knotted} (\textit{i.e.}~knots with high crossing number appear at low grid size) compared to other models. 

\begin{table}
\centering
\begin{tabular}{|c|c|c|c|}
\hline
Knot & probability & occurrences & sample size\\\hline
$5_1$ & 0.006 & 242 & 39200 \\
$7_3$ & 0.000013 & 1 & 76686 \\
$8_8$ & 0.00028 & 13 &  45273 \\
$8_{20}$ & 0.005  & 220 & 41030  \\\hline
\end{tabular}
\caption{Probability of occurences for chirally cosmetic band attachments in our data set~\textit{cf}.~\cite[Table~1]{recent_vaz}}
\label{tab:chiral_probs}
\end{table}

\section{Software and data availability}\label{sec:data_availability}
The code for our computational exploration can be accessed in CK's GitHub repository~\cite{repo}. The software is based on GridPyM~\cite{gridpym}. The repository also contains one example for each band attachment found.

\section{Appendix by Kazuhiro Ichihara, In Dae Jong, Masakazu Teragaito: Chirally cosmetic bandings on knots $7_3$ and $8_8$}

We here give a graphical explanation of why the knots $7_3$ and $8_8$ admit chirally cosmetic bandings. 
Such a banding on $7_3$ is found in this paper and 
that on $8_8$ was discovered in \cite{noteband_vaz} via computer-aided explorations. 
Therefore, it was not immediately clear why such bandings occur.
Also our method can be applied to the banding on the knot $5_1$ given in \cite{zekovic}, although another visual explanation was already provided in \cite{ichiharajongtaniyama}. 

\begin{figure}[htb]
    \centering
    \includegraphics[width=0.8\linewidth]{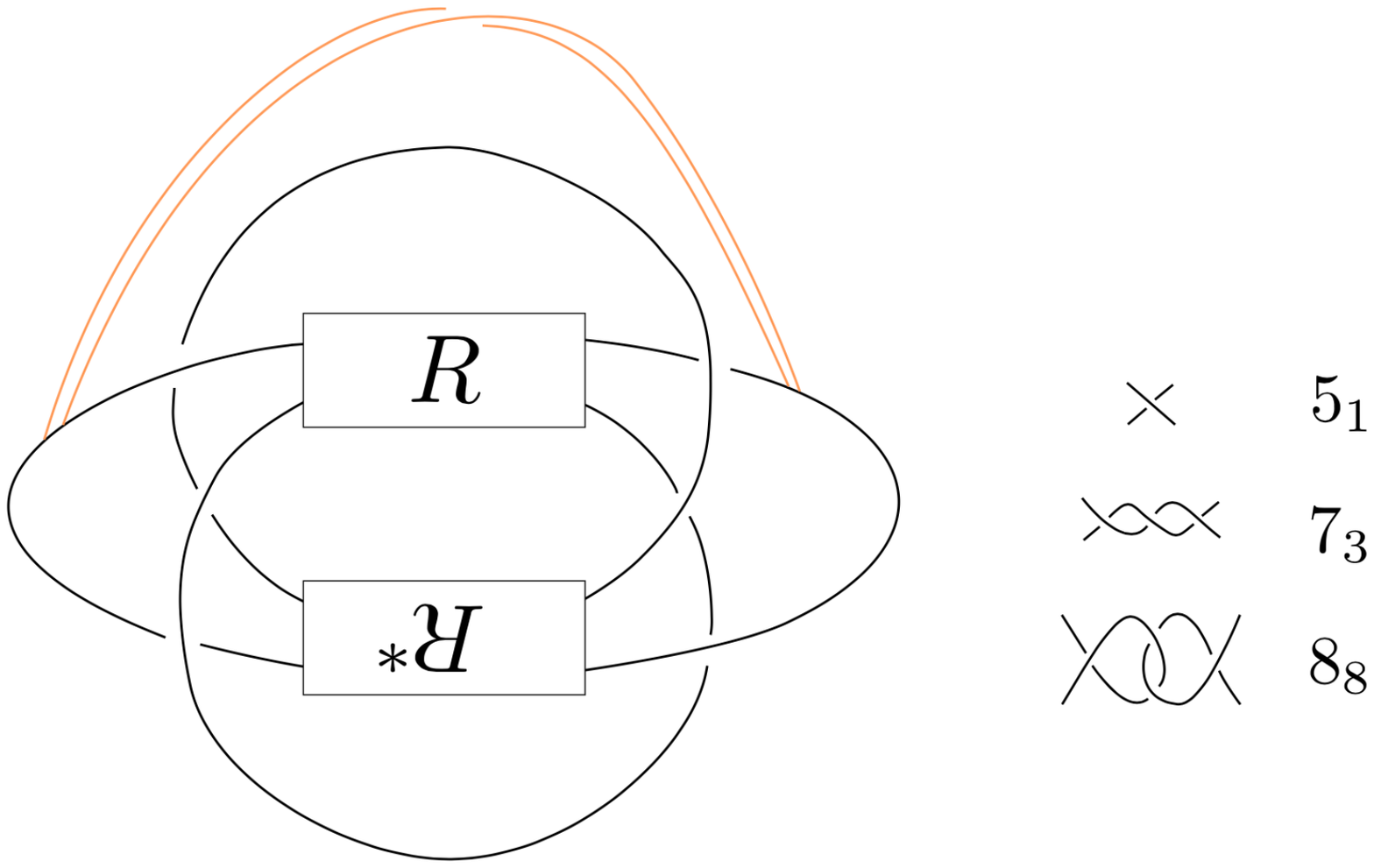}
    \caption{}
    \label{appfig:1}
\end{figure}

We consider a knot diagram as shown in Figure~\ref{appfig:1}. 
The rectangular box $R$ represents a 2-string tangle inserted there, and $R^*$ denotes the mirror image of the tangle $R$. 
We assume that the tangle $R$ admits a reflection symmetry along the vertical axis; that is, $R \cong \reflectbox{$R$}$. 

The diagram of a chirally cosmetic banding on $7_3m$ shown in Figure 5 can be transformed into the diagram in which the box 
$R$ represents a tangle consisting of three half-twists as shown in Figure~\ref{appfig:1}. 
Similarly, the knots $5_1$ and $8_8$ are obtained by replacing $R$ with the respective tangles shown on the right side of Figure~\ref{appfig:1}.

\begin{figure}[htb]
    \centering
    \includegraphics[width=\linewidth]{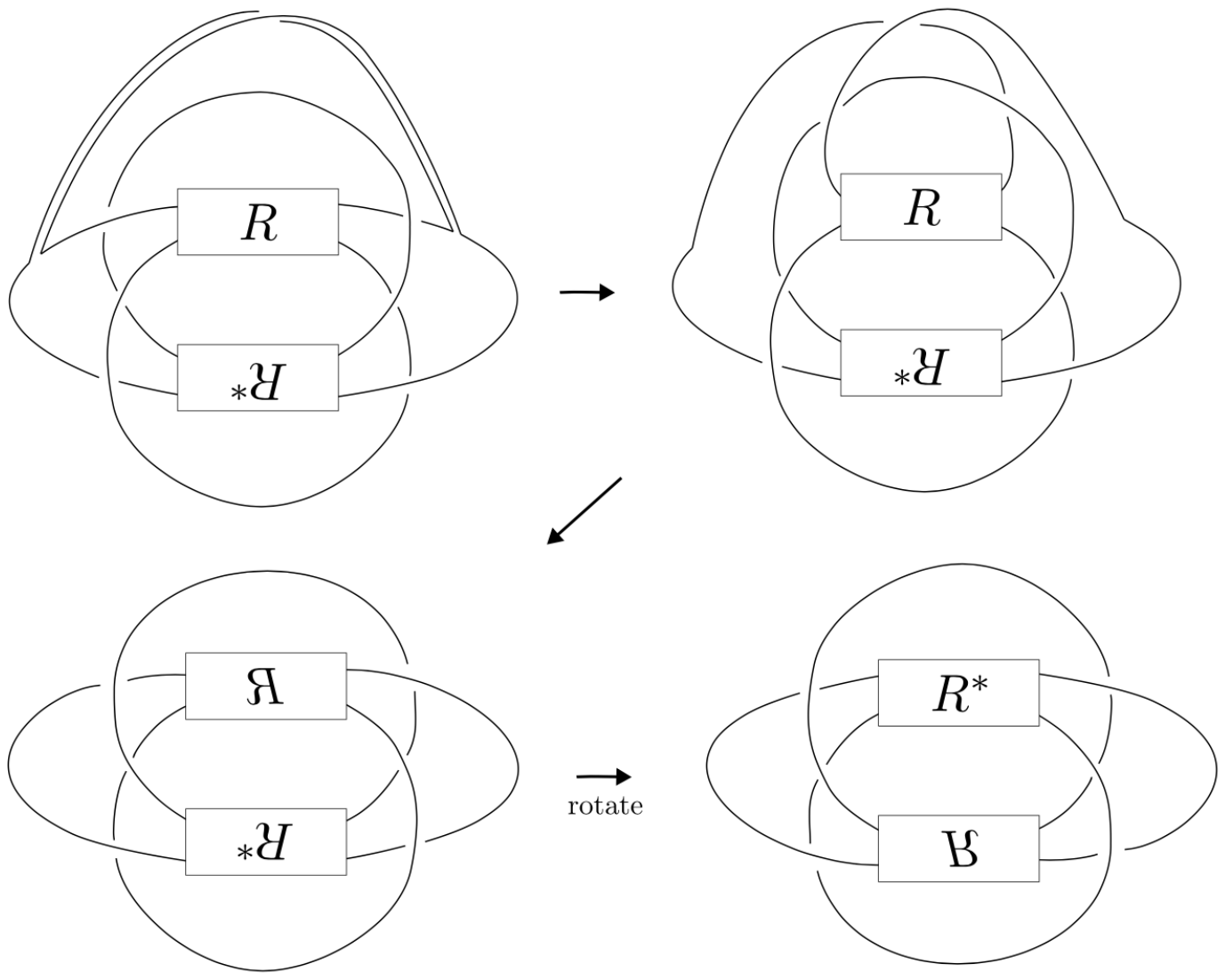}
    \caption{}
    \label{appfig:2}
\end{figure}

Then the banding on the diagram in Figure~\ref{appfig:1} is shown to be chirally cosmetic as illustrated in Figure~\ref{appfig:2}. 
In the figure, moving from the top right to the bottom left, we perform a $\pi$-rotation of only the upper tangle $R$ and the connected arcs along the vertical axis. 
The resultant diagram at the bottom right is actually the mirror image of the original one since we are assuming that $R$ has a reflectional symmetry. 

This trick can give an infinite family of chirally cosmetic bandings on knots.

\newpage
\begin{table}
\centering
\begin{tabular}{|c|p{110mm}|}
\hline
\textbf{knot type} &  \hspace{4cm}{\bf$\hhh$ neighbours}\\\hline\hline
$3_1$ & $0_1$, $3_1$,  $4_1$, $5_2m$, $6_2$, $7_1$, $7_2$, $7_3$, $7_3m$, $7_5$, $7_5m$, $7_6$, $7_7m$, $8_2$, $8_2m$, $8_4m$, $8_6m$, $8_7$, $8_9$, $8_{12}$, $8_{13}$, $8_{13}m$, $8_{14}m$, $8_{15}m$, $8_{19}m$,  $8_{21}m$,  $3_1\#3_1$, $3_1\#3_1m$, $3_1m\#3_1m$, $3_1\#5_1$, $3_1\#5_1m$, $3_1\#5_2$, $3_1\#5_2m$ \\\hline

$4_1$ (A) & $3_1$,  $4_1$,  $5_2$, $6_3$, $7_1$, $7_4$, $7_5$, $8_1$, $8_2$, $8_6$,  $8_7$, $8_{10}$, $8_{11}$, $8_{15}$, $8_{16}$, $8_{17}$, $3_1\#4_1$ \\\hline

$5_1$ & $0_1$, $5_1$, $5_1m$, $6_1m$, $6_2$, $7_2$, $7_2m$, $7_6$, $8_4m$, $8_5$, $8_{14}m$, $8_{20}$, $8_{20}m$, $8_{21}m$, $4_1\#4_1$, $3_1\#5_1$, $3_1m\#5_1$, $3_1m\#5_2m$  \\\hline 

$5_2$ & $0_1$,  $3_1m$, $4_1$, $5_2$, $6_1$, $6_2m$, $6_3$, $7_1m$, $7_5$, $7_5m$, $7_6m$, $7_7$,  $8_1m$, $8_2$, $8_5m$, $8_6$, $8_9$,  $8_{10}$, $8_{11}m$, $8_{14}$, $8_{15}m$, $8_{19}$, $8_{21}$, $3_1m\#5_1m$, $3_1\#5_2$,  $3_1m\#5_2$\\\hline

$6_1$ & $0_1$, $5_1m$, $5_2$, $6_1$, $7_2$, $7_3$, $7_3m$, $7_6m$, $8_7$, $8_8$, $8_8m$, $8_{10}$, $8_{11}m$, $8_{14}$, $8_{20}$, $8_{20}m$, $3_1\#3_1m$\\\hline

$6_2$ & $0_1$, $3_1$, $5_1$, $5_2m$, $6_2$, $6_3$, $7_2$, $7_4m$, $7_5$, $7_5m$, $7_6$, $7_7$, $\color{red} \bf 8_1$, $8_2m$, $8_8$, $8_{10}m$, $8_{12}$, $8_{13}$, $8_{13}m$, $8_{15}m$, $8_{16}m$, $8_{20}$,  $3_1m\#5_2$\\\hline

$6_3$ (A) & $4_1$, $5_2$, $6_2$, $6_3$, $7_1$, $7_6$, $7_7$, $\color{red} \bf 8_1$, $8_4$, $8_{14}$, $\color{red} \bf 8_{15}$, $8_{17}$,  $8_{21}$,   $3_1\#5_1$\\\hline

$7_1$ & $0_1$, $3_1$, $4_1$, $5_2m$, $6_3$, $7_1$, $7_5m$, $7_7m$,  $8_1m$, $8_2$, $8_4m$, $\color{red} \bf 8_{11}$, $\color{red} \bf 8_{15}$\\\hline

$7_2$ & $0_1$, $3_1$, $5_1$, $5_1m$, $6_1$, $6_2$, $7_2$, $8_1$, $8_2m$, $\color{red} \bf 8_6m$, $8_7$, $8_8$, $\color{red} \bf 8_8m$, $\color{red} \bf 8_{11}$, $8_{14}m$, $\color{red} \bf 8_{20}$,  $3_1m\#4_1$ \\\hline

$7_3$ & $0_1$, $3_1$, $3_1m$, $6_1$, $6_1m$, $7_3$, $\color{red} \bf 7_3m$, $8_3$, $8_6m$, $8_7$, $\color{red} \bf 8_7m$, $8_9$, $8_{11}$, $\color{red} \bf 8_{19}$, $\color{red} \bf 8_{19}m$ \\\hline

$7_4$ & $0_1$, $4_1$, $6_2m$, $7_4$, $7_7m$, $8_4$, $8_5m$, $8_{12}$, $8_{13}m$, $8_{16}$, $3_1\#3_1$, $\color{red} \bf 3_1m\#3_1m$\\\hline

$7_5$ & $3_1$, $3_1m$, $4_1$, $5_2$, $5_2m$, $6_2$, $6_2m$, $7_1m$, $7_5$, $\color{red} \bf 8_2$, $8_2m$, $8_6$, $8_{10}$, $\color{red} \bf 8_{10}m$, $8_{11}m$, $8_{12}$, $8_{13}$, $8_{14}$, $8_{17}$\\\hline

$7_6$ & $0_1$, $3_1$, $5_1$, $5_2m$, $6_1m$, $6_2$, $6_3$, $7_6$, $7_7m$, $8_1$, $8_8$, $8_{11}$, $8_{12}$, $8_{13}m$, $8_{14}m$, $8_{15}$, $8_{16}m$, $8_{17}$, $8_{19}m$, $3_1m\#4_1$ \\\hline

$7_7$ &  $3_1m$, $5_2$, $6_2$, $6_3$, $7_1m$, $7_4m$, $7_6m$, $7_7$, $8_4$, $8_{13}$, $8_{14}$, $8_{15}$, $8_{16}m$, $8_{18}$, $3_1\#3_1$\\\hline 
\end{tabular}
\caption{Known $\hhh$-adjacency for prime knots with up to $7$ crossings. The label (A) denotes amphichiral knots. Entries in red represent previously unknown band attachments.}
\label{tab:distance_one1}
\end{table}

\begin{table}
\centering
\begin{tabular}{|c|p{105mm}|}
\hline
\textbf{knot type} &  \hspace{4cm}{\bf$\hhh$ neighbours}\\\hline\hline
$8_1$ & $4_1$, $5_2m$, $\color{red} \bf 6_2$, $\color{red} \bf 6_3$, $7_1m$, $7_2$, $7_6$, $8_1$, $\color{red} \bf 8_{14}m$, $8_{21}m$, $3_1\#5_1m$, $3_1m\#5_2$\\\hline 

$8_2$ & $3_1$, $3_1m$, $4_1$, $5_2$, $6_2m$, $7_1$, $7_2m$, $\color{red} \bf 7_5$, $7_5m$, $8_2$, $8_7m$, $8_{10}$, $8_{13}$, $\color{red} \bf 8_{14}$\\\hline

$8_3$ (A) & $0_1$, $7_3$, $8_3$, $8_{20}$, $3_1\#5_1m$\\\hline

$8_4$ & $0_1$, $3_1m$,  $5_1m$, $6_3$, $7_1m$, $7_4$, $7_7$, $8_4$, $8_{10}$, $8_{13}$, $8_{13}m$, $\color{red} \bf 8_{15}m$,  $3_1m\#5_2$\\\hline

$8_5$ & $0_1$, $5_1$, $5_2m$, $7_4m$, $8_5$, $8_8$, $8_8m$, $\color{red} \bf 8_{19}$, $3_1\#3_1m$, $\color{red} \bf 3_1m\#4_1$\\\hline

$8_6$ & $0_1$, $3_1m$, $4_1$, $5_2$, $\color{red} \bf 7_2m$, $7_3m$, $7_5$, $8_6$, $\color{red} \bf 8_7m$, $8_{10}$, $8_{20}m$, $\color{red} \bf 3_1m\#5_1$, $3_1\#5_2m$\\\hline

$8_7$ & $0_1$, $3_1$, $4_1$, $6_1$, $7_2$, $7_3$, $\color{red} \bf 7_3m$, $8_2m$, $\color{red} \bf 8_6m$, $8_7$, $8_9$, $\color{red} \bf 8_{10}m$, $8_{11}$, $8_{14}m$, $\color{red} \bf 8_{20}$\\\hline

$8_8$ & $0_1$, $6_1$, $6_1m$, $6_2$, $7_2$, $\color{red} \bf 7_2m$, $7_6$, $8_5$, $8_5m$, $8_8$, $8_8m$, $8_9$, $8_{14}m$, $3_1\#5_2$, $4_1\#4_1$\\\hline

$8_9$ (A) & $3_1$, $5_2$, $7_3$, $8_7$, $8_8$, $8_9$, $\color{red} \bf 8_{19}$\\\hline

$8_{10}$ & $0_1$, $4_1$, $5_2$, $6_1$, $6_2m$, $7_5$, $\color{red} \bf  7_5m$, $8_2$,  $8_4$, $8_6$, $\color{red} \bf  8_7m$, $8_{10}$, $8_{11}m$, $8_{12}$, $\color{red} \bf 8_{13}$,  $8_{20}m$ \\\hline

$8_{11}$ & $0_1$, $4_1$, $5_2m$, $6_1m$, $\color{red} \bf 7_1$, $\color{red} \bf 7_2$, $7_3$, $7_5m$, $7_6$, $8_7$, $8_{10}m$, $8_{11}$, $8_{13}m$, $8_{20}$, $\color{red} \bf 3_1\#3_1m$\\\hline

$8_{12}$ (A) & $3_1$, $6_2$,  $7_4$, $7_5$,  $7_6$,  $8_{10}$, $8_{12}$, $3_1\#5_2m$ \\\hline

$8_{13}$ & $3_1$, $3_1m$, $6_2$, $6_2m$, $7_4m$, $7_5$, $7_6m$, $7_7$, $8_2$, $8_4$, $8_4m$, $\color{red} \bf 8_{10}$, $8_{11}m$, $8_{13}$, $8_{14}$, $8_{17}$, $3_1\#4_1$\\\hline

$8_{14}$ & $0_1$, $3_1m$, $5_1m$, $5_2$, $6_1$, $6_3$, $7_2m$, $7_5$, $7_6m$, $7_7$, $\color{red} \bf 8_1m$, $\color{red} \bf 8_2$, $8_7m$, $8_8m$, $8_{13}$, $8_{14}$, $8_{16}$, $\color{red} \bf 8_{19}$, $3_1\#4_1$\\\hline

$8_{15}$ & $3_1m$, $4_1$, $5_2m$, $6_2m$, $\color{red} \bf 6_3$, $\color{red} \bf 7_1$, $7_6$, $7_7$, $\color{red} \bf 8_4m$, $8_{15}$, $8_{21}$, $3_1m\#3_1m$ \\\hline

$8_{16}$ & $0_1$, $4_1$, $6_2m$, $7_4$, $7_6m$, $7_7m$, $8_{14}$, $8_{16}$, $8_{20}m$\\\hline

$8_{17}$ (A) & $4_1$, $6_3$, $7_5$, $7_6$, $8_{13}$,  $8_{17}$, $8_{21}$, $\color{red} \bf 3_1\#5_1$\\\hline

$8_{18}$ (A) & $7_7$, $8_{18}$, $8_{21}$\\\hline

$8_{19}$ & $0_1$, $3_1m$, $5_2$, $\color{red} \bf 7_3$, $\color{red} \bf  7_3m$, $7_6m$, $\color{red} \bf 8_5$, $\color{red} \bf 8_9$, $\color{red} \bf \color{red} \bf 8_{14}$, $8_{19}$, $8_{21}$, $\color{red} \bf 3_1\#3_1m$\\\hline

$8_{20}$ & $0_1$, $5_1$, $5_1m$, $6_1$, $6_1m$, $6_2$, $\color{red} \bf 7_2$, $8_3$, $8_6m$, $\color{red} \bf 8_7$, $8_{10}m$, $8_{11}$, $8_{16}m$, $8_{20}$, $8_{20}m$ \\\hline

$8_{21}$ & $3_1m$, $5_1m$, $5_2$, $6_3$, $8_1m$, $8_{15}$, $8_{17}$, $8_{18}$, $8_{19}$, $8_{21}$\\\hline

\end{tabular}
\caption{Known $\hhh$-adjacency for prime knots with $8$ crossings. The label (A) denotes amphichiral knots. Entries in red represent previously unknown band attachments.}
\label{tab:distance_one2}
\end{table}

\begin{table}
\centering
\begin{tabular}{|c|p{105mm}|}
\hline
\textbf{knot type} &  \hspace{4cm}{\bf$\hhh$ neighbours}\\\hline\hline
$3_1\#3_1$ & $3_1$, $3_1m$, $7_4$, $\color{red} \bf 7_4m$, $ 7_7$, $8_{15}m$, $3_1\#3_1$, $3_1\#4_1$, $3_1\#5_2m$ \\\hline

$3_1\#3_1m$ (A) & $3_1$,  $6_1$, $8_5$, $\color{red} \bf 8_{11}$, $\color{red} \bf 8_{19}$,  $3_1\#3_1m$, $\color{red} \bf 3_1\#4_1$, $\color{red} \bf 3_1\#5_2$\\\hline 

$3_1\#4_1$ & $0_1$, $4_1$, $7_2m$, $7_6m$, $\color{red} \bf 8_5m$, $8_{13}$, $8_{14}$,  $3_1\#3_1$, $\color{red} \bf 3_1\#3_1m$, $4_1\#4_1$, $3_1\#4_1$, $3_1\#5_2$, $3_1\#5_2m$\\\hline

$4_1\#4_1$ (A) &  $5_1$, $8_8$, $3_1\#4_1$, $4_1\#4_1$\\\hline

$3_1\#5_1$ & $3_1$, $5_1$, $5_2m$, $6_3$, $\color{red} \bf 8_{17}$, $3_1\#5_1$, $3_1\#5_1m$\\\hline 

$3_1\#5_1m$ & $3_1$, $5_1m$, $8_1$, $8_3$, $\color{red} \bf 8_6m$, $3_1\#5_1$, $3_1\#5_1m$\\\hline
 
$3_1\#5_2$ & $0_1$, $3_1$, $5_1m$, $5_2$, $8_8$, $\color{red} \bf 3_1\#3_1m$, 
 $3_1\#4_1$, $3_1\#5_2$\\\hline
 
 $3_1\#5_2m$ & $3_1$, $5_2m$, $6_2m$, $8_1m$, $8_4m$, $8_6$, $8_{12}$, $3_1\#3_1$, $3_1\#4_1$, $3_1\#5_2m$  \\\hline 
\end{tabular}
\caption{Known $\hhh$-adjacency for composite knots with up to $8$ crossings. The label (A) denotes amphichiral knots. Entries in red represent previously unknown band attachments.}
\label{tab:distance_one3}
\end{table}

\end{document}